\documentclass{amsart}
\usepackage[utf8]{inputenc}

\usepackage{hyperref}
\usepackage{amsthm}
\usepackage{amssymb,latexsym,amsmath}
\usepackage{enumerate}
\usepackage[all,cmtip]{xy}
\usepackage{graphicx}

\newcommand{\set}{\mathbf{Set}}

\newcommand{\E}{\mathcal{E}}

\DeclareMathOperator{\ub}{ub}
\DeclareMathOperator{\sh}{Sh}
\DeclareMathOperator{\teo}{Th}
\DeclareMathOperator{\ev}{ev}
\DeclareMathOperator{\dom}{dom}
\DeclareMathOperator{\cod}{cod}
\DeclareMathOperator{\Sub}{Sub}

\newtheorem{theorem}{Theorem}[section]
\newtheorem{proposition}[theorem]{Proposition}
\newtheorem{lemma}[theorem]{Lemma}

\newtheorem{axiom}{Axiom}

\theoremstyle{definition}
\newtheorem{definition}[theorem]{Definition}
\newtheorem{remark}[theorem]{Remark}

\newtheorem{example}[theorem]{Example}

\begin{document}

\title{A Categorical Generalization of Counterpoint}




\author[O. A. A. A.]{Octavio A. Agustín-Aquino}
\author[J. S. A. V.]{Juan Sebastián Arias-Valero}
\author[E. R. H.]{Enrique Ruiz Hernández}
\address{Instituto de Física y Matemáticas, Universidad Tecnológica de la Mixteca, Huajuapan de León, Oaxaca, México}
\address{Departamento de Matemáticas, Facultad de Ciencias, Universidad Nacional de Colombia, Bogotá, Colombia}
\address{Centro de Investigación en Teoría de Categorías, México}

\email{octavioalberto@mixteco.utm.mx}
\email{jsariasv1@gmail.com}
\email{e.ruiz-hernandez@cinvcat.org.mx}


\thanks{This work was partially supported by a grant from the \emph{Niels Hendrik Abel Board}.}

\maketitle

\begin{abstract}
We extend Mazzola's counterpoint model using category theory, generalizing from the category $\mathbf{Set}$ to an arbitrary topos other than $\set$. 
This generalization suggests that counterpoint's essential structure depends on specific categorical conditions rather than classical set-theoretic reasoning.

A key contribution is identifying sufficient requirements for a well-behaved counterpoint theory in a topos: some version of Zorn's Lemma (GJZL), and two-valuedness and split supports (NS).

Within a topos, we introduce (weak) quasidichotomies
alongside the classical notion of dichotomy. These structures capture varying degrees of oppositional structure between consonance and dissonance, with weak quasidichotomies preserving the non-Boolean flexibility essential to musical practice while quasidichotomies represent maximal opposition short of complete partition.

We prove a generalized counterpoint theorem giving sufficient conditions for the existence of admitted successors. When the ambient topos turns non-zero successor objects into points, admitted succession can be iterated to form counterpoint paths, which may terminate at consonances with no admitted successor. The framework naturally accommodates counterpoint with sets instead of pure pitches, relaxing the ``yes/no'' character of classical consonance definitions and emphasizing context-dependence.

Mazzola's model allows a Kuratowski closure operator induced by a polarity, which defines an internal topology enabling algebraic-topological analysis of counterpoint structure. We conclude by showing this construction generalizes to involutive morphisms. 
This categorical approach provides foundations for understanding both the historical evolution of contrapuntal practice and cross-cultural divergences in interval organization.
\end{abstract}

\section{Introduction}\label{S:intro}

Our intention in generalizing Mazzola's counterpoint theory \cite{AJM15} via category theory and topos theory is to relativize the notions of consonance and dissonance, and to be able to apply contrapuntal techniques to other musical objects besides pitch.

The nature of the original definitions from Mazzola's counterpoint model is set-theo\-re\-ti\-cal. And rightly so, since Renaissance counterpoint
consolidated the notion of consonance as a Boolean one: an interval is consonant or dissonant, never both or neither. But in the framework of category
theory, Heyting algebras are more natural, and could be used to reflect the historical development of counterpoint via the pseudocomplement. Another approach to non-Booleanness could be by means of fuzzy sets. However,
if we take consonances to define a fuzzy set (as in \cite[Definition 1.2.1]{NW00}), that is, a function $\kappa:S\to [0,1]$, then the pseudocomplement
\[
 \neg \kappa(x) = \begin{cases}
 0, & \kappa(x)>0,\\
 1, & \kappa(x)=0,
 \end{cases}
\]
defines, in particular, a crisp set, even if we apply it to a general fuzzy set. Further applications of this pseudocomplement yield only crisp sets, so we get back to Boolean complements (see \cite[Theorem 2.1.9]{NW00}). This is relevant since the progressive definition of the nature of consonance and dissonance that ended up with a stark separation between them with good contrapuntal properties is a well-known musicological problem (see
\cite{kjS74}); we must point to \cite{gM17} for a discussion of how this could be understood as a mathematical fact discovered by musical means. Hence, the intuitionistic logic of topoi can be used to take the first steps in order to solve this problem.

In particular, the topologization of a counterpoint model via a Kuratowski closure can be put in perspective within the general investigation of closure operators (see \cite{gC03,DT95}, for example). More precisely, Mazzola does all the classical constructions of music in the functor category $[R\text{-}\mathbf{Mod}_{fg},\mathbf{Set}]$ for some ring $R$. We think that a generalization by making $\mathbf{Set}$ any topos satisfying certain conditions and $R$-$\mathbf{Mod}_{fg}$ any (essentially) small category would increase the explanatory power of mathematical musicology. For example, the logic of the internal language of a topos might be more appropriate to reflect the development of counterpoint throughout history via the pseudocomplement. Our definitions are relative to a functor $F:\mathcal{M}\rightarrow\mathcal{E}$ where $\mathcal{M}$ is an (essentially) small category and $\mathcal{E}$ a topos.

The general plan of the article is to construct a categorical generalization of Mazzola's original conception \cite[Part VII]{gM02}, providing examples of certain gains we obtain from it, and thus refining our requirements on the ambient categories required to have a successful generalization; topoi appear as a satisfactory option. We close with a brief study of the Kuratowski operator introduced by Mazzola in order to ``topologize'' counterpoint in the generalized setting. We provide the minimal background of Mazzola's original theory in Section~\ref{S:background}, and we refer the reader to \cite{AJM15}, \cite{kJ92}, and \cite{aM65} for general introductions on these topics.

We generalize the notion of polarity with the concept of quasipolarity in Section~\ref{S:Quasipolarities} in order to introduce the concept of dichotomy; in Subsection~\ref{SS:sdichosnpol}, we introduce the concept of polarity proper, and present the notions of consonance and counterpoint symmetry. In Section~\ref{S:FRondichos}, we present two results for quasipolarities for a functor whose codomain is $\set$. In Section~\ref{S:CPtheorem}, we present a counterpoint theorem. In Section~\ref{S:Conditionsonatopos}, we take stock of all the properties our ambient topos should have in order to have a well-behaved counterpoint theory. Finally in Section~\ref{S:Topology}, we categorize the notion of Kuratowski operator as introduced by Mazzola.

\section{Background}\label{S:background}

Mazzola's counterpoint theory assumes a twelve-tone equal temperament and that pitches are represented up to octave equivalence, and thus we work on the ring $\mathbb{Z}_{12}$. Since the Renaissance, the intervals can be identified with the same ring (for they are all the possible distances modulo octave), and they are split into halves of the same cardinality, the consonances $K=\{0,3,4,7,8,9\}$ and the dissonances $D=K^{c} =\{1,2,5,6,10,11\}$. This stark separation of consonances and dissonances is the final stage of an evolution that is documented in \cite{kjS74}.

Let us describe the affine group 
\[
 \mathrm{Aff}(\mathbb{Z}_{12}):=\mathbb{Z}_{12}\rtimes \mathbb{Z}_{12}^{\times}
\]
as functions that act on $\mathbb{Z}_{12}$, obtained by composing by a multiplication
\[
M_{v}(x) = vx
\]
and a transposition
\[
T^{u}(x) = x+u,
\]
thus any $g\in\mathrm{Aff}(\mathbb{Z}_{12})$ is of the form $T^{u}\circ M_{v}$, with composition given by
\[
 (T^{u_{1}}\circ M_{v_{1}})\circ (T^{u_{2}}\circ M_{v_{2}}) = T^{u_{1}+v_{1}u_{2}}\circ M_{v_{1}v_{2}}.
\]

It turns out that there exists a unique affine automorphism $p=T^{2}\circ M_{5}$ such that
\[
p(K)=D\quad\text{and}\quad p(D)=K.
\]

In particular, $p$ is involutive (meaning $p\circ p = \mathrm{id}$) and a derangement (meaning it has no fixed points).

For the purposes of the most elementary counterpoint, which requires two voices with notes of identical duration, we need to attach an interval $y$ to a cantus firmus $x$. This is done through dual numbers,
\[
\mathbb{Z}_{12}[\epsilon] = \{x+\epsilon y: x,y\in\mathbb{Z}_{12}, \epsilon^{2}=0\}.
\]

This set can also be partitioned into halves of equal cardinality, namely the consonant intervals $K[\epsilon]=\mathbb{Z}_{12}+\epsilon K$ and the dissonant intervals $D[\epsilon]=\mathbb{Z}_{12}+\epsilon D$. For these there exists an affine symmetry $p[\epsilon]=T^{\epsilon 2}\circ M_{5}$ (at cantus firmus $0$) that exchanges consonant and dissonant intervals, and that is involutive.

Thus we can define \emph{counterpoint symmetries} for a consonant interval $\xi\in K[\epsilon]$ as those $g\in \mathrm{Aff}(\mathbb{Z}_{12}[\epsilon])$ that satisfy
\begin{enumerate}
\item $\xi\in g(D[\epsilon])$,
\item $p[\epsilon]\circ g = g\circ p[\epsilon]$,
\item to be the ones where $g(K[\epsilon])\cap K[\epsilon]$ has the largest cardinality among those that satisfy the first two conditions.
\end{enumerate}

Hence $\xi$ is seen as a \emph{deformed dissonant} interval through $g$, and hence it resolves towards a \emph{deformed consonant} interval that is also a consonant interval. The set $g(K[\epsilon])\cap K[\epsilon]$ consists of the \emph{admitted successors}. The Kleiner Kontrapunktsatz \cite[Theorem 2.3]{AJM15} states that every consonant interval has at least $36$ admitted successors.

The maximal intersection condition admits a topological and homological interpretation (see \cite[Chapters 9-10]{AJM15} for further details). The polarity $p[\epsilon]$ induces the Kuratowski closure
\[
\mathrm{cl}(X)=X\cup p[\epsilon](X)
\]
on $\mathbb{Z}_{12}[\epsilon]$. This not only supplies a topology, but also the closure of the set $g(K[\epsilon])\cap K[\epsilon]$ decomposes into orbits of the polarity $\{\eta, p[\epsilon]\eta\}$ where $\eta$ is an admitted successor.

In the counterpoint groupoid these orbits are precisely the curve-connected components of the full subcategory on $\mathrm{cl}(g(K[\epsilon])\cap K[\epsilon])$. Thus we can construct an unnormalized hypergestural homology such that there is one free generator of $H_{1}$, i.e.,
\[
\mathrm{rk}_{R} H_{1}(\mathrm{cl}(g(K[\epsilon])\cap K[\epsilon]);R)=|g(K[\epsilon])\cap K[\epsilon]|.
\]

Our intention in the following is to give a categorical recast of these constructions, so it is possible to understand what happens when the notion of consonance and dissonance is not always an either/or situation (as in Example \ref{Ex:Garlandia}), and also how the dual numbers construction can be either generalized (Example \ref{Ex:SetCounterpoint}) or altogether bypassed (Proposition \ref{T:Counterpoint}).

\section{Quasipolarities and Dichotomies}\label{S:Quasipolarities}

\begin{definition}
Let $\mathcal{E}$ be a nondegenerate topos, $\mathcal{M}$ a small category, and $F:\mathcal{M}\longrightarrow \mathcal{E}$ a functor. A \emph{quasipolarity} for $F$ is a morphism $p:S\longrightarrow S$ of $\mathcal{M}$ satisfying the following conditions:
\begin{enumerate}[i)]
\item The identity $p\circ p=\mathrm{id}_S$ holds.
\item The equalizer $E$ in $\mathcal{E}$ of the pair $F(p),\mathrm{id}_{F(S)}$
\[\xymatrix{
E\ar[r] & F(S)\ar@<.5ex>[r]^{id}\ar@<-.5ex>[r]_{F(p)} & F(S)
}\]
has no points.

If $F(S)=0$, we say $p$ has \emph{trivial image}. 
\end{enumerate}
\end{definition}

\begin{remark}
If the choice of the functor $F$ is obvious (particularly if it is the forgetful functor), then we will omit it.
\end{remark}

This definition was done thinking of the following example.

\begin{example}\label{1}
Take $\mathcal{E}=\textbf{Set}$, $\mathcal{M}=\mathbf{ModAf}_{\mathbb{Z}_{12}}$ (where $\mathbf{ModAf}_R$ denotes the category of modules over a commutative ring with affine transformations between them\footnote{We prefer this notation instead of $\mathbf{Mod}_R$ used by Mazzola, since the latter is the standard notation for the category of $R$-modules and $R$-homomorphisms.}), $F$ the forgetful functor from modules to sets, $S=\mathbb{Z}_{12}$. Then $p=T^{2}\circ M_{5}$
is a quasipolarity.
\end{example}

\begin{example}
The presence of the domain category $\mathcal{M}$ in the definition of quasipolarity is crucial. For instance, if $\mathcal{E}=\mathbf{Set}$, as in the previous example but now 
$\mathcal{M}=\mathbf{FinSet}$, $F$ is the inclusion functor and $S=\{0,1,\ldots,11\}$ (which
coincides with $\mathbb{Z}_{12}$ as a set), then the two permutations
\begin{align*}
 p&=(0,2)(1,7)(3,5)(4,10)(6,8)(9,11),\\
 q&=(0,1)(2,3)(4,5)(6,7)(8,9)(10,11)
\end{align*}
are quasipolarities; $p$ is exactly the same one as the one from the previous example as a function between sets, whereas $q$ does not come from an affine transformation.
\end{example}

\begin{definition}
Let $p:S\longrightarrow S$ be a quasipolarity for a functor $F$. A \emph{dichotomy relative to} $p$ is a pair of monomorphisms
\[\kappa:K\rightarrowtail F(S)\quad\text{and}\quad \delta:D\rightarrowtail F(S)\]
in $\mathcal{E}$ such that
\begin{enumerate}[i)]
\item The canonical morphism from the coproduct $K\amalg D$ (in $\mathcal{E}$) to $F(S)$ is an isomorphism.
\item The monomorphisms $F(p)\circ \kappa$ and $\delta$ are the same subobject of $F(S)$.
\end{enumerate}    
\end{definition}

We denote\footnote{We shall fall into the common abuse of denoting a subobject $A\rightarrowtail B$ by the name of its domain (see before A1.3.1 in \cite{pJ02}, for example).} a dichotomy relative to $p$ by $(\kappa|\delta)_{p}$ or $(K|D)_p$, or simply by $(\kappa|\delta)$ or $(K|D)$ if the election of $p$ is clear.

\begin{example}\label{2}
The subsets
\[
 K=\{0,3,4,7,8,9\}\qquad\text{and}\qquad D=\{1,2,5,6,10,11\}
\]
of $\mathbb{Z}_{12}$ constitute a dichotomy $(K|D)_{T^{2}\circ M_{5}}$ for the polarity of Example \ref{1}.
\end{example}

It is known that the first and foremost consonances were the unison, the fifth and the fourth, while much later the ``imperfect'' thirds and sixths were added. In order to capitalize on the intuitionistic logic of a topos to model the historical development of counterpoint, we will broaden the notion of dichotomy as follows.

\begin{definition}
Let $\E$ be a topos and $P\in\E$ an internal poset of $\E$. A point $p$ of $P$ is said to be a ghostly internal maximal element of $P$ if it satisfies
\[
\forall x\in P.\,p\leq x\Rightarrow\neg\neg(p=x).
\]
Remember that $p$ is an internal maximal element of $P$ if it satisfies $\forall x\in P.\,p\leq x\Rightarrow p=x$.
\end{definition}

\begin{definition}\label{D:weakquasidicho}
A \textit{weak quasidichotomy} relative to a quasipolarity $p:S\rightarrow S$ for a functor $F:\mathcal{M}\rightarrow\mathcal{E}$ is a subobject $K\rightarrowtail FS$ satisfying $Fp(K)\subseteq K^c$, where $K^c$ is the pseudocomplement of $K$ with respect to $FS$. If $K$ is a ghostly internal maximal element in the internal poset $\{u\subseteq FS:Fp(u)\subseteq u^c\}$ of $\E$ then $K$ is said to be a \textit{quasidichotomy}.
\end{definition}

\begin{example}\label{Ex:Garlandia} As reported by Sachs \cite[p. 60]{kjS74}, Johannes de Garlandia (among other 13th-century theorists) considered as consonant (with varying degrees of ``perfection'') what can be represented by the set $K=\{0,3,4,5,7\}$, in particular because he argues that intervals over an octave with respect to a given pitch share the same degree of consonance (or dissonance, for that matter) as those within an octave \cite[Chapter 9]{jG}. The set $K$ does not correspond to a dichotomy, but it is a weak quasidichotomy: for there exist four involutive affine morphisms, namely $T^{1}\circ M_{11}$, $T^{6}\circ M_{1}$, $T^{9}\circ M_{7}$ and $T^{10}\circ M_{5}$, that map $K$ to a subset of the complement $K^{c}=\{1,2,6,8,9,10,11\}$.
\end{example}

\begin{example}\label{Ej:EjemploNoNS} The former example motivates the following. Consider the topos $[C_2\times E,\set]$, where $C_2$ is the cyclic group of order 2 and $E$ is the monoid with exactly two elements: the unit and an idempotent element. The monoid $C_2\times E$ may be described as a commutative monoid with exactly four elements: the unit $e$, an idempotent element $\varphi$, an involutive one $\gamma$ and $\varphi\gamma$. Let $F:C_2\times E\rightarrow\set$ be the following functor. $F(\bullet):=\{0,1,2,3\}=:S$, $F\varphi:=id$ and $F\gamma=(2,3)$, the transposition of 2 and 3. Let $p:S\rightarrow S$ be defined as $p:=(0,1)(2,3)$. Then $p$ is a natural transformation $F\Rightarrow F$.

The equalizer of $p$ and $id_F$ is 0 in $[C_2\times E,\set]$. Let $K:C_2\times E\rightarrow\set$ be the functor $K(\bullet):=\{0\}$ and $Kh:=id$ for every $h\in C_2\times E$. Then $K$ is a subfunctor of $F$ with $K^c(\bullet)=\{1,2,3\}$, $K^c(\varphi)=id$, $K^c(\gamma)=(2,3)$ such that $p(K)$ is a strict subfunctor of $K^c$. There is no greater subfunctor $L$ of $F$ containing $K$ such that $p(L)$ is a subfunctor of $L^c$. This means that $K$ is an external maximal element of $S$.

It is not a ghostly internal maximal element of $S$. To see that, we shall use the approach that a functor $C_2\times E\rightarrow\set$ is viewed as a $C_2\times E$-set. We shall need the following (from now on we shall abbreviate $C_2\times E$ as $C$):
\begin{enumerate}
 \item The subobject classifier of $[C,\set]$ is the $C$-set $\Omega:=\{I\subseteq C\mid I\text{ is a left ideal}\}$; specifically,
 \[
 \Omega=\{\emptyset,\{\varphi,\varphi\gamma\},C\}.
 \]
 Let $\{\varphi,\varphi\gamma\}$ be denoted 1/2. The action of $C$ on $\Omega$ does the following:
 \[
 \zeta\divideontimes I:=\{\xi\in C\mid\xi\zeta\in I\}.
 \]
 For example,
 \[
 e\divideontimes 1/2=1/2,\quad\varphi\divideontimes 1/2=C,\quad\gamma\divideontimes 1/2=1/2\quad\varphi\gamma\divideontimes 1/2=C.
 \]
 The arrow true $t:1\rightarrow\Omega$ is defined as $t(\ast)=C$. The characteristic arrow $\chi_A$ for a sub-$C$-set $A$ of a $C$-set $X$ is defined as
 \begin{align}
 \chi_A(x) :&=\{\zeta\in C\mid\zeta x\in A\}\notag\\
 &=\begin{cases}
     C &\text{ if $x\in A$,}\\
     1/2 &\text{ if $x\notin A$ and $1/2\cdot x\subseteq A$,}\\
     \emptyset &\text{ if $1/2\cdot x\not\subseteq A$.}
 \end{cases}\notag
 \end{align}

 \item An exponential object $Y^X$ is the set $[C,\set](C\times X,Y)$ of equivariant maps with action of $C$ on it defined as $\zeta\odot r$ with
 \[
 \zeta\odot r(\xi,x):=r(\xi\zeta,x)
 \]
 for $r\in [C,\set](C\times X,Y)$. The evaluation map $Y^X\times X\rightarrow Y$ is defined as $\ev_Y^X(r,x):=r(e,x)$. The transpose $\hat{s}$ of an equivariant map $s:X\times Y\rightarrow Z$ is defined by $\hat{s}(x)(\zeta,y):=s(\zeta x,y)$.

 \item For an $X\in[C,\set]$, the arrow $\delta_X:X\times X\rightarrow\Omega$ does the following:
 \[
 \delta_X(\langle x,x'\rangle)=\begin{cases}
  C &\text{ if $x=x'$,}\\
  1/2 &\text{ if $1/2\cdot\langle x,x'\rangle\subseteq \Delta_X$,}\\
  \emptyset &\text{ if $1/2\cdot\langle x,x'\rangle\not\subseteq\Delta_X$,}
 \end{cases}
 \]
 where $\Delta_X=\{\langle x,x\rangle\mid x\in X\}$.

 \item The characteristic arrow $fa:1\rightarrow\Omega$ of $\emptyset\rightarrow 1$ in $[C,\set]$ is defined as $fa(\ast)=\emptyset$. Hence its characteristic arrow $\neg:\Omega\rightarrow\Omega$ does the following:
 \[
 \neg(B)=\begin{cases}
    C &\text{ if $B=\emptyset$,}\\
    \emptyset &\text{ if $B=1/2,C$.}
 \end{cases}
 \]

 \item Given a $C$-set $X$, the characteristic arrow of $1_X$, denoted $t_X:X\rightarrow\Omega$, is the constant function $(x\mapsto C)$. The characteristic $\forall_X:\Omega^X\rightarrow\Omega$ of `$t_X$', the name of $t_X$, is as follows:
 \[
 \forall_X(r)=\begin{cases}
    C &\text{ if $r=t_{C\times X}$,}\\
    1/2 &\text{ if $r\neq t_{C\times X}$ and $1/2\odot r=t_{C\times X}$,}\\
    \emptyset &\text{ if $1/2\odot r\neq t_{C\times X}$.}
 \end{cases}
 \]

 \item The arrow $\wedge:\Omega\times\Omega\rightarrow\Omega$ is the expected one. The same holds for the equalizer of $\wedge$ and $\pi_1$. So the order on $\Omega$ is
 \[
 \emptyset\leq 1/2\leq C.
 \]

 \item The arrow $\Rightarrow:\Omega\times\Omega\rightarrow\Omega$ is the function defined as
 \[
 A\Rightarrow B=\begin{cases}
 C &\text{ if $A\leq B$,}\\
 1/2 &\text{ if $A=C$ and $B=1/2$,}\\
 \emptyset &\text{ if $A=C,1/2$ and $B=\emptyset$.}
 \end{cases}
 \]

 \item Since $[C,\set](C\times X,\Omega)\cong\Sub(C\times X)$, we can think of $\Omega^X$ as the set of sub-$C$-sets of $C\times X$.

 \item The formula $x\in u\Rightarrow x\in v$ in the internal logic is interpreted in the topos as the arrow
 \[\xymatrix@C=.5cm{
 \Omega^X\times\Omega^X\times X\ar[r]_(.45){1\times\Delta_X} & \Omega^X\times\Omega^X\times X\times X\ar[rr]_{1\times twist\times 1} & & \Omega^X\times X\times\Omega^X\times X\ar[r]_(.67){\ev_X\times\ev_X} & \Omega\times\Omega\ar[r]_(.6)\Rightarrow & \Omega.
 }\]
 Denote that arrow by $\psi$. The formula $u\subseteq v$ in the internal logic is interpreted in the topos as the arrow
 \[\xymatrix{
 \Omega^X\times\Omega^X\ar[r]^(.6){\hat{\psi}} & \Omega^X\ar[r]^{\forall_X} & \Omega.
 }\]
 Define $\preccurlyeq:=\forall_X\circ\hat{\psi}$. From the previous item in the list, we get
 \[
 A\preccurlyeq B=\begin{cases}
    C &\text{ if $A\subseteq B$,}\\
    1/2 &\text{ if $A\not\subseteq B$ and $1/2\cdot A\subseteq B$,}\\
    \emptyset &\text{ if $1/2\cdot A\not\subseteq B$.}
 \end{cases}
 \]

 \item The arrow $\Omega^p:\Omega^X\rightarrow\Omega^X$ is the transpose of $\ev\circ(1\times p)$. By calculating it we get $\Omega^p(r)=r\circ(1\times p)$; that is, given $A\subseteq C\times S$, we have $\Omega^p(A)=(1\times p)^{-1}(A)$. Hence its left adjoint is the direct image under $1\times p$.

 \item A sub-$C$-set $A$ of a $C$-set $X$ corresponds to a point of $\Omega^X$, viz. the sub-$C$-set $C\times A$ of $C\times X$, and vice versa. Said differently: the points of $\Omega^X$ are the fixed points of $\Omega^X$ under the action of $C$.
\end{enumerate}
Now, given the item 8 in the list, the order relation $R$ on $\Omega^X\times\Omega^X$, calculated as the pullback of $\preccurlyeq$ along $t:1\rightarrow\Omega$, may be thought of as an order relation on $\Sub(C\times X)\times\Sub(C\times X)$ induced by the isomorphism in the item 8.

Let $D$ be the object $\{u\subseteq S:p(u)\subseteq u^c\}$ in $[C,\set]$. It may be thought of as a subset of $\Sub(C\times S)$. Hence the internal order in $D$ induced by that of $\Omega^S$ may be thought of as its external order induced by that of $\Sub(C\times S)$.

The sub-$C$-set $\{0\}$ of $S$ corresponds to the sub-$C$-set $C\times\{0\}$ of $C\times S$, and $C\times\{0\}$ is not maximal in $D$ since
\[
(C\times\{0\}\preccurlyeq C\times\{0\}\cup\{\langle e,2\rangle,\langle\varphi,2\rangle,\langle\varphi\gamma, 3\rangle,\langle\gamma,3\rangle\})=C
\]
and
\[
\neg(C\times\{0\}=C\times\{0\}\cup\{\langle e,2\rangle,\langle\varphi,2\rangle,\langle\varphi\gamma, 3\rangle,\langle\gamma,3\rangle\})=C.
\]
So put
\[
I:=C\times\{0\}\cup\{\langle e,2\rangle,\langle\varphi,2\rangle,\langle\varphi\gamma, 3\rangle,\langle\gamma,3\rangle\}
\]
and
\[
J:=C\times\{0\}\cup\{\langle \varphi\gamma,2\rangle,\langle\gamma,2\rangle,\langle e, 3\rangle,\langle\varphi,3\rangle\}.
\]
Let $E:=\{I,J\}$. Thus $E$ is a sub-$C$-set of $D$ since
\[
e\odot\chi_I=\chi_I,\quad \varphi\odot\chi_I=\chi_I,\quad\varphi\gamma\odot\chi_I=\chi_J,\quad\gamma\odot\chi_I=\chi_J 
\]
and
\[
e\odot\chi_J=\chi_J,\quad \varphi\odot\chi_J=\chi_J,\quad\varphi\gamma\odot\chi_J=\chi_I,\quad\gamma\odot\chi_J=\chi_I. 
\]
It is clear that $E\rightarrow 1$ is epic.

Now, with abuse of notation, let $K:1\rightarrow D$ be the point $C\times\{0\}$ of $D$, and let $j:E\hookrightarrow D$ be the inclusion of $E$ into $D$. Therefore if $L:=\langle K\circ !_E,j\rangle:E\rightarrow D\times D$, then the arrow
\[\xymatrix{
E\ar[r]_(.38)L & D\times D\ar@{^{(}->}[r] & \Sub(C\times S)\times\Sub(C\times S)\ar[r]_(.68)\cong & \Omega^S\times\Omega^S\ar[r]_{\langle\preccurlyeq,\neg\circ\delta_{\Omega^S}\rangle} & \Omega\times\Omega\ar[r]_(.6)\wedge & \Omega
}\]
is the constant arrow sending everything to $C$. Therefore, by Kripke-Joyal semantics,
\[
\exists u\in D.\{0\}\subseteq u\wedge\neg(\{0\}=u)
\]
is true in the internal logic of the topos.
\end{example}

\begin{example}
Consider the Sierpiński topos $\set^{\bullet\rightarrow\bullet}$, which is localic. Let $j:4\rightarrow 6$ be the inclusion function of $4=\{0,1,2,3\}$ into $6=\{0,\dots,5\}$. Now, let $p:4\rightarrow 4$ be the following composite of transpositions: $(0,2)(1,3)$, and let $q:6\rightarrow 6$ be the same composite $(0,2)(1,3)$. The following diagram commutes:
\[\xymatrix{
4\ar[d]_j\ar[r]^p & 4\ar[d]^j \\
6\ar[r]_q & 6.
}\]
The equalizer of $id_j$ and $\langle p,q\rangle$ is the function $0\rightarrow\{4,5\}$, which has no points.

Now, let $k=id_2$. Hence $\langle p,q\rangle(k)\subseteq k^c$, where $k^c$ is the inclusion of $\{2,3\}$ into $\{2,3,4,5\}$. It is clear that $k$ is a maximal subfunction  of $j$ such that $\langle p,q\rangle(k)\subseteq k^c$, and that $\langle p,q\rangle(k)\neq k^c$. That is, $k$ is an external maximal element of the object $\{u\subseteq j:p(u)\subseteq u^c\}$. Let us show it is an internal maximal element. For that we need the following:

\begin{enumerate}
 \item The subobject classifier of $\set^{\bullet\rightarrow\bullet}$ is the function $\Omega:3\rightarrow 2$ satisfying
\[\Omega(n):=\begin{cases}
0 &\text{if $n=0$,}\\
1 &\text{if $n\neq 0$,}
\end{cases}\]
where 3 is thought of as $\{0,1/2,1\}$. The arrow true $t:1\rightarrow\Omega$ is the pair $\langle t^u,t^d\rangle$ of functions $t^u:1\rightarrow 3$ and $t^d:1\rightarrow 2$ defined as $t^u(\ast)=1$ and $t^d(\ast)=1$.

 \item An exponential object $h^g$ in $\set^{\bullet\rightarrow\bullet}$ is the function
\[\xymatrix{
\set^{\bullet\rightarrow\bullet}(g,h)\ar[r]^(.45){h^g} & \set(\cod g,\cod h)
}\]
\[\xymatrix{
\langle r,s\rangle\ar@{}[r]|(.6){\mapsto} & s;
}\qquad\quad\]
the evaluation map $\ev:h^g\times g\rightarrow h$ is the pair of functions $\langle\ev^u,\ev^d\rangle$ defined as follows:
\[\ev^u(\langle r,s\rangle, c):=rc\quad\text{and}\quad\ev^d(w,d)=wd,\]
with $\ev^u:\set^{\bullet\rightarrow\bullet}(g,h)\times\dom g\rightarrow\dom h$ and $\ev^d:\set(\cod g,\cod h)\times\cod g\rightarrow\cod h$, and the transpose of $\langle u,v\rangle:f\times g\rightarrow h$ is defined as the pair
\[\langle\langle u^t,v^t\circ f\rangle,v^t\rangle:f\rightarrow h^g,\]
where $u^t$ and $v^t$ are the transposes of $u$ and $v$ in $\set$.

\item Given $f\in\set^{\bullet\rightarrow\bullet}$, $\delta_f:f\times f\rightarrow\Omega$ is a pair of functions $\delta_f^u:A\times A\rightarrow 3$, $\delta_f^d:B\times B\rightarrow 2$ defined as follows:
\[\delta_f^u(a,a')=\begin{cases}
1 &\text{if $a=a'$,}\\
1/2 &\text{if $a\neq a'$ and $fa=fa'$,}\\
0 &\text{if $fa\neq fa'$;}
\end{cases}\]
and
\[\delta_f^d(b,b')=\begin{cases}
1 &\text{if $b=b'$,}\\
0 &\text{if $b\neq b'$.}
\end{cases}\]

 \item Given a function $f:A\rightarrow B$, the characteristic of $1:f\rightarrow f$, denoted $t_f:f\rightarrow\Omega$, is the pair $\langle t_f^u,t_f^d\rangle$, where $t_f^u:A\rightarrow 3$ and $t_f^d:B\rightarrow 2$ are the constant functions $(a\mapsto 1)$ and $(b\mapsto 1)$, respectively. The characteristic $\forall_f:\Omega^f\rightarrow\Omega$ of `$t_f$', the name of $t_f$, is as follows:
\[\forall_f^u(\langle r,s\rangle)=\begin{cases}
1 &\text{if $r=t_f^u$ and $s=t_f^d$,}\\
1/2 &\text{if $r\neq t_f^u$ and $s=t_f^d$,}\\
0 &\text{if $s\neq t_f^d$,}
\end{cases}\]
and
\[\forall_f^d(h)=\begin{cases}
1 &\text{if $h=t_f^d$,}\\
0 &\text{if $h\neq t_f^d$.}
\end{cases}\]

 \item The arrow $\wedge:\Omega\times\Omega\rightarrow\Omega$ is the pair $\langle\wedge^u,\wedge^d\rangle$ of functions $\wedge^u:3\times 3\rightarrow 3$ and $\wedge^d:2\times 2\rightarrow 2$ with the expected result as if 3 was ordered as $0\leq 1/2\leq 1$ and 2 as $0\leq 1$.

 The equalizer
 \[\xymatrix{
 \leq_3\ar@{^{(}->}[r]\ar[d]_\leq & 3\times 3\ar@<.5ex>[r]^(.6){\wedge^u}\ar@<-.5ex>[r]_(.6){\pi_1}\ar[d]_{\Omega\times\Omega} & 3\ar[d]^\Omega \\
 \leq_2\ar@{^{(}->}[r] & 2\times 2\ar@<.5ex>[r]^(.6){\wedge^d}\ar@<-.5ex>[r]_(.6){\pi_1} & 2
 }\]
 in $\set^{\bullet\rightarrow\bullet}$ yields the order expressed above for 3 and 2.
 \item The arrow $\Rightarrow:\Omega\times\Omega\to\Omega$ is the pair of functions $\langle\Rightarrow^u,\Rightarrow^d\rangle$ where $\Rightarrow^d$ is the usual material implication, and $\Rightarrow^u$ is defined as follows:
 \[
 n\Rightarrow^u m=\begin{cases}
 1&\text{ if $n\leq m$,}\\
 1/2&\text{ if $n=1$ and $m=1/2$,}\\
 0&\text{ if $n=1/2,1$ and $m=0$.}
 \end{cases}
 \]

 \item Since $\dom(\Omega^f)=\set^{\bullet\rightarrow\bullet}(f,\Omega)\cong\Sub(f)$ and $\Omega^f(\langle r,s\rangle)=s$, we can think of $\Omega^f$ as the set of subfunctions $\Sub(f)$ of $f$.

 \item The formula $x\in u\Rightarrow x\in v$ in the internal logic is interpreted in the topos as the arrow
 \[\xymatrix@C=.8cm{
 \Omega^f\times\Omega^f\times f\ar[r]_(.45){1\times\Delta_f} & \Omega^f\times\Omega^f\times f\times f\ar[rr]_{1\times twist\times 1} & & \Omega^f\times f\times\Omega^f\times f\ar[r]_(.67){\ev_f\times\ev_f} & \Omega\times\Omega\ar[r]_(.6)\Rightarrow & \Omega.
 }\]
 Denote that arrow by $\psi$. The formula $u\subseteq v$ in the internal logic is interpreted as the arrow
 \[\xymatrix{
 \Omega^f\times\Omega^f\ar[r]^(.6){\hat{\psi}} & \Omega^f\ar[r]^{\forall_f} & \Omega.
 }\]
 Define $\subseteq:=\forall_f\circ\hat{\psi}$. It is the pair $\langle\subseteq^u,\subseteq^d\rangle$ such that, given the previous item in the list,
 \[
 g\subseteq^ug'=\begin{cases}
 1 &\text{ if $A\subseteq A'$ and $B\subseteq B'$}\\
 1/2 &\text{ if $A\not\subseteq A'$ and $B\subseteq B'$}\\
 0 &\text{ if $B\not\subseteq B'$,}
 \end{cases}
 \]
 where $g:A\rightarrow B$ and $g':A'\rightarrow B'$ are subfunctions of $f$.
 It is clear what $\subseteq^d$ should be.
\end {enumerate}
Now, given the item 7 in the list above, the order relation $R$ on $\Omega^f\times\Omega^f$, calculated as the pullback of $\subseteq$ along $t:1\rightarrow\Omega$, may be thought of as an order relation on $\Sub(f)\times\Sub(f)$ induced by the isomorphism in the item 7.

Let $A$ be the object $\{u\subseteq j:\langle p,q\rangle(u)\subseteq u^c\}$ in $\set^{\bullet\rightarrow\bullet}$. It may be thought of as a subset of $\Sub(j)$, where $\langle p,q\rangle(g)$ is calculated pointwise for every $g\in\Sub(j)$. Hence $A\times A$ is a subset of $\Sub(j)\times\Sub(j)$. That is, the internal order in $A$ induced by that of $\Omega^j$ may be thought of as its external order induced by that of $\Sub(j)$.

We claim that the arrow
\[\xymatrix{
A\ar[r]^(.38){\langle k!,1\rangle} & A\times A\ar@{^{(}->}[r] & \Omega^j\times\Omega^j\ar[rr]^{\langle\subseteq,\delta_{\Omega^j}\rangle} & &\Omega\times\Omega\ar[r]^(.6)\Rightarrow & \Omega
}\]
is the constant 1, where $\subseteq$ is the arrow defined in item 8. Indeed, the only way this arrow can fail to be constant 1 is if the consequent of the resultant formula is 0, and the antecedent is $1/2$ or 1 for some subfunction $g$ of $j$ in $A$, or if the consequent is $1/2$ and the antecedent is 1 for some subfunction $g$ of $j$ in $A$.

Now, let $\langle r_k,s_k\rangle$ and $\langle r_g,s_g\rangle$ be the characteristic arrows of $k$ and some subfunction $g$ of $j$ in $A$, respectively, and notice that there is no function in $A$ above $k$. So if $\delta_{\Omega^j}^u(k,g)=0$, then $s_k\neq s_g$; if $k\subseteq^u g=1$, then $\langle r_k,s_k\rangle=\langle r_g,s_g\rangle$, which is a contradiction; if $k\subseteq^u g=1/2$, then, as there is no function in $A$ above $k$, we would have $\cod k=\cod g$, and accordingly $s_k=s_g$, a contradiction again. Therefore it is not possible to have the first two cases. If $\delta_{\Omega^j}^u(k,g)=1/2$, then $r_k\neq r_g$ and $s_k=s_g$; if $k\subseteq^u g=1$ then $r_k=r_g$, a contradiction.

Therefore the formula $k\subseteq u\Rightarrow k=u$ for $u$ a variable of type $A$ is true in the internal logic of the topos. Hence so is $\forall u\in A.k\subseteq u\Rightarrow k=u$ (see the discussion in the first lines in \cite[\S 18.2]{cM01}).
\end{example}

\subsection{Strong Dichotomies and Polarities}\label{SS:sdichosnpol}

\begin{definition}
Let $(\kappa|\delta)$ be a dichotomy relative to a quasipolarity $p: S\longrightarrow S$. We say that the dichotomy is \emph{strong} if $p$ is the unique quasipolarity such that $F(p)\circ \kappa$ and $\delta$ are the same subobject; in that case, we say that $p$ is a \emph{polarity}.
\end{definition}

\begin{example}
The dichotomy of Example \ref{2} is strong and its quasipolarity is a polarity. \qed
\end{example}

\begin{example}
The notion of polarity also stresses the importance of the choice of the category $\mathcal{M}$ and the functor\footnote{A purely combinatorial
relativization was explored in \cite{oA12}.} $F:\mathcal{M}\rightarrow\mathcal{E}$. For instance, if $\mathcal{M}_{1}=\mathbf{ModSAf}_{\mathbb{Z}_{12}}$, which is the
subcategory of modules over a commutative ring with the morphisms restricted to those of the form $T^{u}\circ M_{\pm 1}$, taking
\[
(\kappa|\delta)=(\{0,2,3,4,7,8\}\mid\{1,5,6,9,10,11\}),
\]
we have that $p=T^{1}\circ M_{-1}$ is a polarity, but in the supercategory $\mathcal{M}_{2}=\mathbf{ModAf}_{\mathbb{Z}_{12}}$ it is not since
\[
 T^{1}\circ M_{-1}(\kappa)=T^{9}\circ M_{7}(\kappa)=\delta.
\]
\end{example}

  
\begin{definition} Let $(\kappa|\delta)$ be a dichotomy. A \emph{consonance} is a morphism $\xi$ (of $\mathcal{E}$) from the terminal object $\mathbf{1}$ of $\mathcal{E}$ to $K$, that is, a point of $K$. More generally, a \emph{generalized consonance} can be defined as a morphism $\xi$ of $\mathcal{E}$ with codomain $K$. 
\end{definition}

\begin{definition}
Let $(\kappa|\delta)$ be a dichotomy with quasipolarity $p:S\longrightarrow S$. A \emph{counterpoint symmetry} for a consonance $\xi:\mathbf{1}\longrightarrow K$ is an isomorphism $g:S\longrightarrow S$ of $\mathcal{M}$ such that
\begin{enumerate}[i)]
\item the morphism \label{D1} $\kappa\circ \xi$ is a point of
$F(g)\circ \delta$, or, more precisely, $\kappa\circ \xi$ factors through $F(g)\circ \delta$.
\item \label{D2} the identity $g\circ p=p\circ g$ holds, and
\item \label{D3} if an isomorphism $g':S\longrightarrow S$ of $\mathcal{M}$ satisfies \ref{D1} and \ref{D2}, then there is a monomorphism in $\E$
from 
$(F(g')\circ \kappa)\cap\kappa$ to $(F(g)\circ \kappa)\cap\kappa$.
\end{enumerate}
\end{definition}

\begin{remark}
The subobjects $F(g)\circ \kappa$ and $F(g)\circ \delta$ of $F(S)$ represent the
\emph{$g$-deformed consonances} and \emph{dissonances}, respectively.
\end{remark}

\begin{remark}
Since $\kappa\circ\xi$ is a deformed dissonance, then the tension induced by $g$ has to be relaxed towards a consonance that is also a $g$-deformed consonance. Thus the subobject $(F(g)\circ \kappa)\cap\kappa$ corresponds to the \emph{admitted successors}. More precisely, a consonance $\eta:\mathbf{1}\longrightarrow K$ is an admitted successor of the consonance $\xi$ if the morphism $\kappa\circ\eta$ factorizes through $(F(g)\circ \kappa) \cap \kappa$.
\end{remark}

In the case when $\mathcal{E}=\mathbf{Set}$ and $\mathcal{M}=\mathbf{ModAf}_{\mathbb{Z}_{12}[\epsilon]}$,  the condition \ref{D2} is exactly the same requirement stated
by the classical model of Mazzola. The condition \ref{D1} captures the idea of revealing the tension
between consonances and dissonances by making $\xi$
a $g$-deformed dissonance, while condition \ref{D3} implies a maximum of artistic choices.

\begin{remark}
If $p$ is a polarity, then it provides a canonical candidate for a counterpoint symmetry for every consonance $\xi:\mathbf{1}\to K$: condition \ref{D1} follows from $F(p)(K)=D$, while condition \ref{D2} is immediate. If $p$ is the only isomorphism satisfying these two conditions, then it is automatically a counterpoint symmetry. In this case, however,
\[
(F(p)\circ\kappa)\cap\kappa=\delta\cap\kappa=0,
\]
so there are no admitted successors. Thus the categorical model naturally allows \emph{cul-de-sacs}: consonances for which a counterpoint symmetry exists but its admitted-successor object is initial.

Thus the class of isomorphisms satisfying conditions \ref{D1} and \ref{D2} is non-empty whenever a polarity exists. In general, however, we do not claim that this class necessarily contains an element satisfying the maximality condition \ref{D3}; whenever such an element exists, it is a counterpoint symmetry in the sense defined above.

On the other hand, whenever a counterpoint symmetry $g$ for $\xi$ has a non-zero admitted-successor object
\[
(F(g)\circ\kappa)\cap\kappa,
\]
and this object has a point $\eta:\mathbf{1}\to K$, the construction may be repeated with $\eta$ in place of $\xi$. Hence admitted successors may be iterated until, possibly, a \textit{cul-de-sac} is reached.
\end{remark}

More generally, a \emph{counterpoint symmetry} for a generalized consonance $\xi:E\longrightarrow K$ is an isomorphism $g:S\longrightarrow S$ of $\mathcal{M}$ such that  $\kappa\circ Im(\xi)$ factors through $F(g)\circ \delta$, plus the conditions \ref{D2} and \ref{D3} above. In this case, we require the existence of the image $Im(\xi)\rightarrowtail K$ of $\xi$.

\begin{remark}
The notion of \emph{cantus firmus} and \emph{discantus} is obtained for general $R$-modules as follows. Let $R$ and $S$ be commutative rings. The restriction $r:R\to S$ gives rise to the functor 
\begin{align*}
 S\otimes_{R}-:\mathbf{ModAf}_{R}&\to\mathbf{ModAf}_{S}\\
 N&\mapsto S\otimes_{R} N,\\
 f=T^{k}\circ f_{0}:N\to K &\mapsto S\otimes_{R}f=T^{1\otimes
 k}\circ S\otimes_{R}f_{0}.
\end{align*}

Recall now that the \emph{dual numbers} for a ring
$R$ is the quotient
\begin{align*}
 R[\epsilon]:=\frac{R[x]}{\langle x^{2}\rangle}
 = \{a+\epsilon b:a,b\in R\text{ and } \epsilon^{2}=0\}.
\end{align*}

We have the restriction $i:R\to R[\epsilon]:a \mapsto a+\epsilon 0$, which allows us to
construct, for an $R$-module $M$, the module $M[\epsilon]=R[\epsilon]\otimes_{R}M$, which is
the \emph{dual numbers module with respect to
$M$}. As $R$-modules, we have
\[
M[\epsilon] \cong M\oplus M,
\]
and thus the first component of $(m_{1},m_{2})\in M[\epsilon]$ corresponds to the cantus firmus and $m_{2}$ is the \emph{interval} that separates it from the discantus. The composition of this functor with the forgetful functor $F:\mathbf{ModAf}_{\mathbb{Z}_{12}[\epsilon]}\to \mathbf{Set}$ recovers the classical Mazzola model for counterpoint via counterpoint symmetries.
\end{remark}

\begin{example}\label{Ex:SetCounterpoint}
Let $S\neq \emptyset$ be a finite set. We know that $\wp(S)$ can
be seen as a $\mathbb{Z}_{2}$-algebra (the so-called \emph{Boolean algebra}) defining the product by intersection and addition by symmetric difference
\[
 A+_{\wp(S)}B := A\triangle B = (A\cup B)\setminus (A\cap B).
\]

In particular,
\[
 p = T^{S}\circ M_{1}
\]
is a quasipolarity that coincides with the operation of calculating the complement in $S$. For a given family $\kappa$ of $2^{|S|}/2=2^{|S|-1}$ of subsets of $S$, we define
\[
\mathrm{Stab}(\kappa)=\{U\subseteq S:T^{U}(\kappa)=\kappa\},
\]
so we have that $T^{S}\circ M_{1}$ is a polarity if, and only if,
$\mathrm{Stab}(\kappa)=\{\emptyset\}$. Denote by $\delta=\complement_{\wp(S)} \kappa$ the set of dissonances. Now we can see that $p$ is polarity whenever $\mathrm{Stab}(\kappa)=\{\emptyset\}$: suppose there is a set $U$ such that for any $K\in \kappa$
\[
T^{U}(K) = U\triangle K = D,
\]
for some $D\in \delta$. Since $U\triangle K = S\triangle K'$ for some $K'\in \kappa$, then $(S\triangle U) \triangle K = K'$, but by hypothesis
$S\triangle U = \emptyset$, which means $U=S$.

If $p=T^{S}\circ M_{1}$ is a polarity and $U$ is a subset of $S$ such that $T^{U}(\kappa)=\kappa$, then
\[
 T^S\circ T^U(\kappa) = T^{S\triangle U}(\kappa) = T^S(T^U(\kappa))=T^S(\kappa)=\delta
\]
and uniqueness of $p$ implies $S\triangle U = S$. Hence $U=\emptyset$, which means $\mathrm{Stab}(\kappa)=\{\emptyset\}$.

Non-empty dichotomies with $\mathrm{Stab}(\kappa)=\{\emptyset\}$ exist. Suppose $S$ is a set of cardinality $2k+1$ for $k\in\mathbb{N}$. Then take
\[
\kappa = \{U\subseteq S:|U|\leq k\}.
\]

If for any $K\in \kappa$ we have
\[
 T^{U}(K)=U\triangle K\in \kappa
\]
then $T^{U}(\emptyset) = U\in\kappa$, for $\emptyset\in \kappa$. Now $|S\setminus U|\geq k+1$, which means we can choose $A\subseteq S\setminus U$ with $|A|=k$, i.e., $A\in \kappa$. Now
\[
 |T^{U}(A)|=|A\triangle U| = |A|+|U| \geq k
\]
and since $T^{U}(A)\in \kappa$, it means $|U|=0$, namely $U=\emptyset$.

Thus we can construct the algebra $\wp(S)[\epsilon]$ such that the cantus firmus and the interval are sets. The sets of consonances and dissonances in $\wp(S)[\epsilon]$ are, of course,
\[
 \kappa[\epsilon]=\wp(S)+\epsilon \kappa\quad\text{and}\quad \delta[\epsilon]=\wp(S)+\epsilon \delta.
\]

As in the case of the classical counterpoint model, it can be proved that the calculation of
counterpoint symmetries can be reduced to cantus firmus $0$ and the quasipolarity
\[
p_{0}=T^{\epsilon S}\circ M_{1}.
\]

The affine morphism $g=T^{\epsilon  U}\circ M_{1+\epsilon W}$ always commutes with $p_{0}$, and for a consonance $\xi = \epsilon K\in \kappa[\epsilon]$ it will occur that it is a $g$-deformed dissonance whenever there is $C+\epsilon D\in \delta[\epsilon]$ such that
\[
 \xi=\epsilon K=g(C+\epsilon D)=C+\epsilon (U\triangle D\triangle (W\cap C)),
\]
thus $C=0$ and
\begin{equation}\label{Eq:Hichert}
 K = U\triangle D.
\end{equation}

From this point on it is straightforward to generalize Hichert's algorithm \cite[Algorithm 2.1]{AJM15} to complete the calculations, for we only have to iterate over the candidates $g=T^{\epsilon  U}\circ M_{1+\epsilon W}$ that satisfy \eqref{Eq:Hichert}, determine the cardinality of $g(\kappa[\epsilon])\cap \kappa[\epsilon]$ and keep the ones that maximize the intersection.
\end{example}

\section{Further Remarks on Dichotomies}\label{S:FRondichos}

\begin{proposition}
If $F:\mathcal{M}\rightarrow\mathbf{Set}$ is a functor with $\mathcal{M}$ a small category, then every quasipolarity $p:S\rightarrow S$ for $F$ has a dichotomy associated to it.
\end{proposition}
\begin{proof}
Let $p:S\rightarrow S$ be a quasipolarity. Consider $Fp:FS\rightarrow FS$. If $FS$ is finite, it is easy to see that $p$ has a dichotomy associated to it.

Suppose $FS$ is infinite. Consider the following relation over $FS$: given $s,t\in FS$, we say $s\sim t$ if and only if $s=Fp(t)$. Let $\sim'$ be the equivalence relation generated by $\sim$. Now, by the Axiom of Choice, form $K$ by choosing one element from each equivalence class under the equivalence relation $\sim'$. Let $D:=Fp(K)$. Then $Fp(D)=K$ and $FS=K+D$.
\end{proof}

Now, given the notion of quasidichotomy, one could ask if in $\set$ a quasidichotomy automatically becomes a dichotomy. The answer is yes.

\begin{proposition}\label{P:quasidichoisdicho}
Let $F:\mathcal{M}\rightarrow\set$ be a functor with $\mathcal{M}$ a small category. Let $p:S\rightarrow S$ be a quasipolarity for $F$ with no trivial image. If $K\subseteq FS$ is a quasidichotomy relative to $p$, then $(K|Fp(K))$ is a dichotomy.
\end{proposition}
\begin{proof}
Suppose $K$ is a quasidichotomy. Then $Fp(K)\subseteq K^c$ and $K$ is not contained in a larger set with the same property. Let us suppose $Fp(K)$ is not $K^c$. Then there exists $x\in K^c$ that is not in $Fp(K)$. Hence
\[
Fp(K)\subseteq K^c\setminus\{x\}=K^c\cap\{x\}^c=(K\cup\{x\})^c.
\]
Now, consider $K\cup\{x\}$. It cannot be the case that
$Fp(K\cup\{x\})=Fp(K)\cup\{Fp(x)\}$ is contained in $(K\cup\{x\})^c$, since $K$ is maximal, thus $Fp(x)\in K\cup\{x\}$, but $Fp$ has no fixed points. Hence $Fp(x)\in K$. Therefore $F(pp)x=x$ is in $Fp(K)$, which is contradictory. Hence $Fp(K)=K^c$ for any maximal $K$.
\end{proof}


\section{A Counterpoint Theorem}\label{S:CPtheorem}

With these preliminaries, we can obtain a rudimentary counterpoint theorem, i.e., a sufficient criterion for a counterpoint symmetry to admit successors.

\begin{lemma}
Let $\mathcal{E}$ be a topos, and $K_1,K_2$ be subobjects of $S\in\mathcal{E}$. If $K_2$ is a complemented subobject of $S$ and $K_1$ is not a subobject of $K_2$ then $K_1\cap K_2^c\neq 0$, where $K_2^c$ denotes the complement of $K_2$.
\end{lemma}
\begin{proof}
If $K_1\cap K_2^c=0$, then, since $K_2$ is complemented, $K_1$ is a subobject of $K_2$, which is a contradiction.
\end{proof}

The following gives a sufficient condition for a counterpoint symmetry not to lead to a \textit{cul-de-sac}: this condition gives us the chance to have at least one counterpoint symmetry not mapping all deformed consonances into dissonances.

\begin{proposition}[Counterpoint theorem]\label{T:Counterpoint}
Let $(K|D)$ be a dichotomy with quasipolarity $p:S\rightarrow S$, and $g$ be a counterpoint symmetry for some consonance. If $F(g)(K)$, the image of $K$ under $Fg$, is not a subobject of $D$, then $F(g)(K)\cap K\neq 0$.
\end{proposition}
\begin{proof}
Notice that $D$ is a complemented subobject of $FS$ with $D^c=K$. Hence, by the previous lemma, $F(g)(K)\cap K\neq 0$.
\end{proof}

\begin{remark}
If the functor $F$ has $\mathbf{Set}$ as codomain and $F(S)$ is finite, then if $F(g)K$ is different from $D$ (a weaker condition) we still have $F(g)K\cap K\neq\emptyset$.
\end{remark}

\begin{remark}
It is useful to regard admitted succession as defining paths between consonances. A finite \emph{counterpoint path} is a sequence of consonances
\[
\xi_0,\xi_1,\ldots,\xi_n
\]
such that, for each $0\leq i<n$, there is a counterpoint symmetry $g_i$ for $\xi_i$ for which $\xi_{i+1}$ is an admitted successor. Once a point of the admitted-successor object has been chosen, it is again a consonance and the construction can therefore be iterated. A path terminates when the admitted-successor object associated with the current counterpoint symmetry is $0$; we also call such a termination a \emph{cul-de-sac}.

In particular, in a topos satisfying NS, the conclusion
\[
F(g_i)(K)\cap K\neq 0
\]
of Proposition~\ref{T:Counterpoint} provides an actual point and hence permits the path to be continued by one further step.
\end{remark}

\section{Conditions on the Ambient Topos}\label{S:Conditionsonatopos}

\begin{definition}
Let $\E$ be a topos and $A\in\E$. Given a term $w$ of power type $PA$ in the internal language of $\E$, define
\[
\ub(w):=\{x:\forall y\in w\,(y\leq x)\},
\]
i.e. $\ub(w)$ is the setlike term of upper bounds for $w$.
Given a poset $(A,\leq)$ in $\E$, $A$ is said to be \textit{inductive} in $\E$ if it satisfies
\[
\forall w\in PA.(w\text{ is a chain})\Rightarrow\exists x\in\ub(w).
\]
\end{definition}

The idea behind Definition~\ref{D:weakquasidicho} is to construct the object $A$ determined by the formula
\[
\{u\subseteq FS:Fp(u)\subseteq u^c\},
\]
where $p$ is a quasipolarity and $Fp(u)$ is the image of $u$ under $Fp$, and to find a ghostly internal maximal object in $A$ since $A$ is ordered by inclusion. Now, this object $A$ is an inductive poset in $\E$ (see Lemma~\ref{L:Aisinductive} below). So we are interested in topoi $\E$ having the following property which we call \textit{Ghostly Johnstone-Zorn's Lemma}:

\begin{axiom}[GJZL]
If $(A,\leq)$ is an inductive poset in $\E$, then there exists $a:1\rightarrow A$ which is a ghostly internal maximal element of $A$.
\end{axiom}

This axiom is inspired by the one appearing in \cite[p. 996]{pJ02}, which is:

\begin{axiom}[JZL]
If $(A,\leq)$ is an inductive poset in $\E$, then there exists $a:1\rightarrow A$ which is an internal maximal element of $A$.
\end{axiom}

The discussion in the previous definition raises the following question: are there topoi other than $\set$ satisfying GJZL? The answer is yes; it is established in the following proposition.

\begin{proposition}\label{P:zornslemmaholdsinlocalics}
Assume Zorn's Lemma holds in $\set$. If $\E$ is a localic $\set$-topos, then $\E$ satisfies JZL (\cite[D4.5.14]{pJ02}), and therefore GJZL too.
\end{proposition}

Proposition~\ref{T:Counterpoint} gives a sufficient condition for the existence of admitted successors of consonances in a general ambient topos. However, this does not guarantee that $F(g)K\cap K$ has points. If we want that, the following property for a topos $\E$ could be desirable:

\begin{axiom}[NS]\label{A:NS}
For every object $A\in\E$, either $A\cong 0$ or $A$ has a point. Equivalently, the topos is two-valued and supports split therein.
\end{axiom}

Since we know that localic $\set$-topoi satisfy JZL, we might want to know if there are non-degenerate localic $\set$-topoi other than $\set$ satisfying NS. The answer is no:

\begin{proposition}
No localic $\set$-topos satisfies NS but $\set$ and the degenerate topos.
\end{proposition}
\begin{proof}
Let $\sh(H)$ be a localic $\set$-topos (i.e. $H$ is a locale). 
Every representable functor in $[\mathcal{O}(H)^{op},\set]$ is a sheaf. But, by Yoneda, every representable $\mathcal{O}(H)(-,U)$ has no points except $\mathcal{O}(H)(-,1)$.
\end{proof}

Given that non-degenerate localic $\set$-topoi other than $\set$ satisfy JZL but not NS, it is natural to ask if there are topoi other than $\set$ satisfying GJZL and NS? The answer is yes. Before the promised result, we need the following lemma:

\begin{lemma}
Let $\E$ be a topos satisfying NS and $P$ an internal poset in $\E$. Then every point $p$ of $P$ is external maximal iff it is a ghostly internal maximal. 
\end{lemma}
\begin{proof}
Suppose $p$ is an external maximal. Consider the object $B:=\{x\in P:p < x\}$, where the expression `$y<x$' means `$y\leq x\wedge\neg(y=x)$'. Then since $\E$ satisfies NS, $B=0$ or $B\neq 0$. Suppose $B\neq 0$. Then $B$ has a point $b$, but then $p$ is not an external maximal element of $P$, which is a contradiction. Thus we get $\neg(B\neq 0)$. Thus from $B=0$ or $B\neq 0$, and $\neg(B\neq 0)$, it follows $B=0$. That is, in the internal language we have
\[
x\in P, p\leq x,\neg(p=x)\vdash_{\teo(\E)}\bot.
\]
Hence $x\in P, p\leq x\vdash_{\teo(\E)}\neg\neg(p=x)$. That is,
\[
\vdash_{\teo(\E)}\forall x\in P\,.p\leq x\Rightarrow\neg\neg(p=x).
\]

Conversely, suppose $p$ is a ghostly internal maximal element of $P$. So let $q$ be another point of $P$ such that $p\leq q$. Then we get the closed term $\neg\neg(p=q)$ of type $\Omega$. But the deductive system of the internal logic of $\E$ is complete (see Theorem 4.31 in \cite{jB08}). So we have $p=q$ or $\neg(p=q)$. Therefore $p=q$ (cf. the discussion in Remark~\ref{R:discussioderived}).
\end{proof}

\begin{theorem}\label{T:GrottoposprecohoversetsatsZL}
Let $p:\mathcal{E}\rightarrow\mathbf{Set}$ be the canonical geometric morphism for a Grothendieck topos $\mathcal{E}$. Assume Zorn's Lemma holds in $\mathbf{Set}$. If $p$ is precohesive over $\set$, then $\E$ satisfies GJZL and NS.
\end{theorem}
\begin{proof}
$\E$ satisfies NS  (see Motivation in \cite{RS25}).

Now, consider the external poset $C$ of chains in $A$, ordered by inclusion. As $\E$ has arbitrary unions of subobjects since it is defined over $\set$ and $C$ is an inductive poset in $\set$, $C$ has an external maximal element $A_0$. Hence
\[
A_0\in PA\wedge A_0\in C.
\] 
Whence, since $A$ is inductive, we have $\exists x\in\ub(A_0)$. That is, $\ub(A_0)\rightarrow 1$ is epic. Hence $\ub(A_0)\neq 0$. So, since $\E$ satisfies NS, $\ub(A_0)$ has a point $a:1\rightarrow\ub(A_0)$. It is clear that $A_0\cup\{a\}$ is a chain with $A_0\subseteq A_0\cup\{a\}$, but $A_0$ is maximal. Therefore
\[
A_0\cup\{a\}=A_0.
\]
This means $a\in A_0$. Whence $a$ is an external maximal element of $A$. By the previous lemma, $a$ is a ghostly internal maximal element of $A$.
\end{proof}

\begin{example}
Some examples of non-Boolean presheaf topoi precohesive over $\set$ are
\begin{itemize}
\item $\set^E$, where $E$ is the monoid with exactly two elements: the unit and an idempotent element;
\item the topos $\widehat{\Delta_1}$ of reflexive graphs;
\item $[(\mathbf{FinOrd}-0)^{op},\set]$, where $\mathbf{FinOrd}$ is the category of finite ordinals whose arrows are all functions between them.
\end{itemize}

More examples can be found in \cite[4.2 and 4.3]{mM14} and in \cite[7.1 and 8.5]{mM21}. Some examples of non-Boolean sheaf topoi precohesive over $\set$ can be found in \cite[1.7 and 12.4]{mM21}.
\end{example}

So far we have not proved that $\{u\subseteq FS:Fp(u)\subseteq u^c\}$ is an inductive internal poset.

\begin{lemma}\label{L:Aisinductive}
Let $p$ be a quasipolarity for a functor $F:\mathcal{M}\rightarrow\E$. Let
\[
A:=\{u\subseteq FS:Fp(u)\subseteq u^c\}.
\]
Then $A$ is an inductive poset in $\E$.
\end{lemma}
\begin{proof}
We proceed à la Bell (see \cite{jB08}). In the internal language we have:
\begin{equation}\label{D:Fpreservesunions}
w\subseteq A\vdash_{\teo(\E)}Fp(\cup w)=\cup_{u\in w} Fp(u).
\end{equation}

Now, we want to show
\begin{equation}\label{D:imagecontained}
w\text{ is a chain of $A$}\vdash_{\teo(\E)}Fp(\cup w)\subseteq(\cup w)^c.
\end{equation}
It is clear that for different variables $u,v$ of the appropriate type
\begin{align}\label{D:totality}
w\text{ is a chain of $A$},u&\in w,v\in w,y\in Fp(u)\notag\\&\vdash_{\teo(\E)}w\subseteq A\wedge u\in w\wedge v\in w\wedge(u\subseteq v\vee v\subseteq u)\wedge y\in Fp(u).
\end{align}
On the other hand,
\begin{align}
w\subseteq A,u\in w,v\in w, u\subseteq v, y\in Fp(u) &\vdash_{\teo(\E)}y\in Fp(u)\subseteq Fp(v)\subseteq v^c\notag\\
&\vdash_{\teo(\E)}y\in v^c,\notag
\end{align}
and
\begin{align}
w\subseteq A,u\in w,v\in w,v\subseteq u,y\in Fp(u) &\vdash_{\teo(\E)}u^c\subseteq v^c\wedge y\in Fp(u)\subseteq u^c\notag\\
&\vdash_{\teo(\E)}y\in v^c.\notag
\end{align}
Thus, by combining the last two derivations, we get
\begin{equation}\label{D:yinvc}
w\subseteq A\wedge u\in w\wedge v\in w\wedge(u\subseteq v\vee v\subseteq u)\wedge y\in Fp(u)\vdash_{\teo(\E)}y\in v^c.
\end{equation}
Therefore, by concatenating Derivations~\eqref{D:totality} and \eqref{D:yinvc}, by Cut, we get
\[
w\text{ is a chain of $A$},u\in w,v\in w,y\in Fp(u)\vdash_{\teo(\E)}y\in v^c.
\]
Hence, by 3.4.2 in \cite{jB08},
\[
w\text{ is a chain of $A$},u\in w,y\in Fp(u)\vdash_{\teo(\E)}\forall v\in w.y\in v^c.
\]
Whence, by 3.7.2 in \cite{jB08},
\[
w\text{ is a chain of $A$},\exists u\in w.y\in Fp(u)\vdash_{\teo(\E)}\forall v\in w.y\in v^c;
\]
i.e., we get Derivation~\eqref{D:imagecontained}.

Now, by Derivations~\eqref{D:Fpreservesunions} and \eqref{D:imagecontained} internally we have
\begin{align}
w\in PA\text{ is a chain} &\vdash_{\teo(\E)}\cup w\in A\wedge\cup w\in\ub(w)\notag \\
&\vdash_{\teo(\E)}\cup w\in\ub(w),
\end{align}
that is, we have
\[
w\in PA\text{ is a chain}\vdash_{\teo(\E)}\exists u\subseteq FS\, (u\in\ub(w)).
\]
Hence
\[
\vdash_{\teo(\E)}\forall w\in PA.(w\text{ is a chain})\Rightarrow\exists u\subseteq FS\,(u\in\ub(w)).
\]
So, by \cite[Proposition 3.30]{jB08}, $A$ is an inductive poset in $\E$.
\end{proof}

\begin{remark}\label{R:discussioderived}
If the codomain category in Proposition~\ref{P:quasidichoisdicho} is a topos other than $\set$, then the conclusion is just $K^c\setminus Fp(K)=0$. So, given the fact that when working in $\set$, it is difficult (if not impossible) to distinguish between internal  and external arguments, there are two possibilities for the proof of the previous proposition: it relies on the principle of excluded middle to arrive at the conclusion $K^c\setminus Fp(K)=0$, i.e. it relies on the truth of `$K^c\setminus Fp(K)=0$ or $\neg(K^c\setminus Fp(K)=0)$' (equivalently, on the truth of `$K^c\setminus Fp(K)=0$  follows from $\neg\neg(K^c\setminus Fp(K)=0)$'); or it relies on the strong-witnessedness of the internal logic of $\set$ since $\set$ is two-valued and supports split in it (see Theorem 4.31 and Proposition 4.32 in \cite{jB08}, and Axiom~\ref{A:NS} NS above).

We elaborate a little bit more on the second possibility. Suppose the functor $F$ has as codomain a non degenerate possibly non-Boolean topos $\E$ satisfying NS. Then for $K^c\setminus Fp(K)$ we would have `$K^c\setminus Fp(K)=0$' or `$K^c\setminus Fp(K)$ has a point', i.e. `$K^c\setminus Fp(K)=0$' or `$\neg(K^c\setminus Fp(K)=0)$'; differently said using the turnstile of the deductive system of the internal logic of the topos:
\begin{equation}\label{D:completeness}
\vdash_{\teo(\E)}K^c\setminus Fp(K)=0\qquad\text{or}\qquad\vdash_{\teo(\E)}\neg(K^c\setminus Fp(K)=0)
\end{equation}
(this is also true because the system is complete, see \cite[Proposition 4.32]{jB08}). We suppose $\vdash_{\teo(\E)}\neg(K^c\setminus Fp(K)=0)$ and we get a contradiction, but the system is consistent since the topos is nondegenerate. Hence we get
\begin{equation}\label{D:consistence}
\not\vdash_{\teo(\E)}\neg(K^c\setminus Fp(K)=0).
\end{equation}
So, from \eqref{D:completeness} and \eqref{D:consistence}, we get $\vdash_{\teo(\E)}K^c\setminus Fp(K)=0$. We did not use the principle of the excluded middle, either outside or inside the system, but NS. (Notice that $K^c\cap\{x\}^c=(K\cup\{x\})^c$ is valid in intuitionistic logic). So we have:
\end{remark}

\begin{proposition}\label{P:JZLnNSimpliesqdisd}
If $\E$ is a nondegenerate topos satisfying NS, then every quasidichotomy relative to a quasipolarity $p$ for a functor $F:\mathcal{M}\rightarrow\E$ is such that $K^c\setminus Fp(K)=0$.    
\end{proposition}

We notice that, even if the topos in Example~\ref{Ej:EjemploNoNS} contains an object that could be part of a quasidichotomy, the topos is not localic and does not satisfy NS. Indeed: 

\begin{example}
Consider Example~\ref{Ej:EjemploNoNS}. It is not hard to see that the subobject classifier $\Omega:C_2\times E\rightarrow\set$ of $[C_2\times E,\set]$ is the following:
\[\Omega(\bullet)=\{0,\textstyle{\frac{1}{2}},1\},\qquad\Omega(e)=\Omega(g)=id,\qquad\Omega(f)(x)=\Omega(fg)(x)=\begin{cases}
 0 &\text{if $x=0$,}\\
 1 &\text{if $x=\frac{1}{2},1$,}
\end{cases}\]
and $\top:1\Rightarrow\Omega$ the natural transformation choosing 1 in $\Omega(\bullet)$. The functor $\mathbf{2}$, which is $\{0,1\}$ on $\bullet$ and the identity arrow on every arrow $h$ in $C_2\times E$ is a retract of $\Omega$. Therefore, by D4.6.2 in \cite{pJ02}, $[C_2\times E,\set]$ is a (non-Boolean) De Morgan topos.

The subterminal objects of $[C_2\times E,\set]$ are 0 and 1. If we consider the functor $G:C_2\times E\rightarrow\set$ defined as $G(\bullet)=\{0,1,2,3\}$, $Gf=id$ and $Gg=(2,3)$, and the natural transformation $q:G\Rightarrow G$ defined as $q(\bullet)=(2,3)$, then the point $1\Rightarrow G$ choosing 0 in $G(\bullet)$ equalizes $q$ and $id$ but $q\neq id$. Therefore, by C1.4.7 and before A1.2.5 in \cite{pJ02}, $[C_2\times E,\set]$ is not localic.
Notice that $\set^E$ is De Morgan (the proof of this is similar to that for the De-Morgan-ness of $[C_2\times E,\set]$).

The topos $[C_2\times E,\set]$ does not satisfy Axiom NS. Indeed, notice that the functor $H:C_2\times E\rightarrow\set$ defined as $H(\bullet):=\{0,1\}$, and $H(f):=id$ and $H(g):=(0,1)$ has no points: it is not possible to have an arrow $r:1\Rightarrow H$ such that
\[\xymatrix{
1(\bullet)\ar[r]^{r(\bullet)}\ar[d]_1 & H(\bullet)\ar[d]^{Hg}\\
1(\bullet)\ar[r]_{r(\bullet)} & H(\bullet)
}\]
commutes.
\end{example}

\section{Topology and closure operators}\label{S:Topology}

Given the results in the previous section, in this one we shall work in a Grothendieck topos\footnote{This implies that, for each morphism of $\mathcal{E}$, the associated direct image functor is left adjoint to the associated inverse image functor, which in turn implies that the direct image functor preserves joins of subobjects (see \cite[Lemma A.1.3.1]{pJ02}).} $\E$ whose canonical geometric morphism is precohesive.

\subsection{The operator induced by an involutive automorphism}\label{Mazzolaop}
The Kuratowski operator was introduced by Mazzola in counterpoint \cite[Section 10.2.2]{AJM15} to topologize counterpoint intervals, so that homology could be calculated and related to the conditions determining the rules of first-species counterpoint. This construction can be generalized with no difficulty to the categorical case whenever we have an involutive morphism. Let $f:E\longrightarrow E$ be a morphism of $\mathcal{E}$ such that $f\circ f=\mathrm{id}_{E}$.
Define an operator on subobjects $M$ of $E$ by means of the equation
\[
 \overline{M}:=M\cup f(M).
\]
We check that this is indeed a Kuratowski closure.

\begin{enumerate}
\item $\overline{0}=0\cup f(0)=0$
\item If $m:M\rightarrowtail E$ is a subobject of $E$ in  $\mathcal{E}$, then $M\subseteq M\cup f(M)= \overline{M}$.
\item If $m:M\rightarrowtail E$ is a subobject of $E$ in  $\mathcal{E}$, then
\begin{align*}
 \overline{\overline{M}}&=\overline{M}\cup f(\overline{M})\\
 &=M\cup f(M)\cup f(M\cup f(M))\\
 &=M\cup f(M)\cup f\circ f(M)\\
 &=M\cup f(M)=\overline{M}.
\end{align*}
\item If $m:M\rightarrowtail E$ and $m':M'\rightarrowtail E$ are subobjects of $E$ in  $\mathcal{E}$, then
\begin{align*}
 \overline{M\cup M'}&=M\cup M'\cup f(M\cup M')\\
 &=M\cup M'\cup f(M)\cup f(M')=\overline{M}\cup \overline{M'}.
\end{align*}
\end{enumerate}

\subsection{The closure operator on subobjects induced by an arbitrary endomorphism}

As we have mentioned regarding the historical development of the notion of consonance and dissonance leading to Renaissance counterpoint (see \cite[p. 60]{kjS74}), not always all intervals could be placed neatly into either category, but some approximate morphism is available. We analyze this case in this section.

Let $f:E\longrightarrow E$ be an arbitrary endomorphism of $\mathcal{E}$. Define an operator on subobjects $M$ of $E$ by means of the equation $\overline{M}:=M\cup f(M)$, where $f(M)$ denotes the image of $M$ under $f$. This operator satisfies the following properties:
\begin{enumerate}[i)]
\item If $m:M\rightarrowtail E$ is a subobject of $E$ in  $\mathcal{E}$, then
\[
 M\subseteq M\cup f(M)= \overline{M}.
\]

\item If $m:M\rightarrowtail E$ and $m':M'\rightarrowtail E$ are subobjects of $E$ in  $\mathcal{E}$, then
\begin{align*}
 \overline{M\cup M'}&=M\cup M'\cup f(M\cup M')\\
 &=M\cup M'\cup f(M)\cup f(M')\\
 &=\overline{M}\cup \overline{M'}.
\end{align*}

\end{enumerate} 

In particular, this operator preserves the order on subobjects, though it needs not be idempotent\footnote{See \cite[p. xiv]{DT95} for another general discussion.}. However, by iterating this operator, we can obtain an idempotent operator satisfying the previous properties, that is, a Kuratowski operator. In fact, if the lattice of subobjects of $E$ is complete (for example, if $\mathcal{E}$ is small-cocomplete), then the operator defined by the equation
\[\overline{M}=\bigcup_{k\in \mathbb{N}}f^{k}(M)=M\cup f(M)\cup f^2(M)\cup \cdots\]
satisfies the following properties:
\begin{enumerate}[i)]
\item Since the direct image along $f$ preserves the least subobject,
\[
\overline{0}=0\cup \bigcup_{k\in \mathbb{N}}f^{k}(0)=0.
\]
\item If $m:M\rightarrowtail E$ is a subobject of $E$ in  $\mathcal{E}$, then $M\subseteq \overline{M}$.
\item If $m:M\rightarrowtail E$, then 
\begin{align*}
\overline{\overline{M}}&=\bigcup_{k\in \mathbb{N}}f^{k}(\overline{M})\\
&=\bigcup_{k\in \mathbb{N}}f^{k}(\bigcup_{j\in \mathbb{N}}f^{j}(M))\\
&=\bigcup_{k\in \mathbb{N}}\bigcup_{j\in \mathbb{N}}f^{k+j}(M)\\
&=\bigcup_{k\in \mathbb{N}}f^{k}(M)=\overline{M}.
\end{align*}
\item If $m:M\rightarrowtail E$ and $m':M'\rightarrowtail E$ are subobjects of $E$ in  $\mathcal{E}$, then 
\begin{align*}
\overline{M\cup M'}&=\bigcup_{k\in \mathbb{N}}f^{k}(M\cup M')\\
&=\bigcup_{k\in \mathbb{N}}f^{k}(M)\cup \bigcup_{k\in \mathbb{N}}f^{k}(M')\\
&=\overline{M}\cup\overline{M'}.
\end{align*}
\end{enumerate} 

In the case when $f^{n}=\mathrm{id}_{E}$ occurs for
some $n\in\mathbb{N}$, this collapses to
\[\overline{M}=M\cup f(M)\cup f^2(M)\cup \cdots \cup f^{n-1}(M)\]
and, of course, if $n=2$, then we obtain Mazzola's closure operator from the previous subsection. Furthermore, given an endomorphism $p:S\longrightarrow S$ of $\mathcal{M}$, the definition above applies to the endomorphism $F(p):F(S)\longrightarrow F(S)$, and hence we have a Kuratowski operator on the subobjects of $F(S)$. 

\section{Final remark}\label{S:Final}

A natural question for future work is: if $S$ is a $K$-finite object of a topos $\E$, $p:S\rightarrow S$ is an arrow in $\E$ such that $pp=1$ and the equalizer of $p$ and $1_S$ is 0, then is it true that for any subobject $A\rightarrowtail S$ maximal such that $pA$ is a subobject of $A^c$, we have $pA=A^c$?

\section*{Acknowledgment}

We thank Guerino Mazzola for his valuable feedback on early versions of this paper.


\end{document}